\documentclass{article}
\usepackage{amsmath}
\usepackage{amssymb}

\begin{document}

\begin{center}
{\Large An efficient method for solving equations in generalized quaternion
and octonion algebras }%
\begin{equation*}
\end{equation*}

\textbf{Cristina Flaut and Vitalii Shpakivskyi}%
\begin{equation*}
\end{equation*}
\end{center}

\textbf{Abstract.} {\small Quaternions often appear in wide areas of applied
science and engineering such as wireless communications systems, mechanics,
etc. It is known that are two types of non-isomorphic generalized quaternion
algebras, namely: the algebra of quaternions and the algebra of
coquaternions. In this paper, we present the formulae to pass from a basis
in the generalized quaternion algebras to a basis in the division
quaternions algebra or to a basis in the coquaternions algebra and vice
versa. The same result was obtained for the generalized octonion algebra.
Moreover, we emphasize the applications of these results to the algebraic
equations and De Moivre's formula in generalized quaternion algebras and in
generalized octonion division algebras.\medskip }


\begin{equation*}
\end{equation*}

\vspace{3mm} \textbf{0. Introduction}\vspace{3mm}%
\begin{equation*}
\end{equation*}

Let $\gamma _{1},\gamma _{2}\in \mathbb{R}\setminus \{0\},$ let $\mathbb{H}%
\left( \gamma _{1},\gamma _{2}\right) $ be the generalized quaternion{\small %
\ }algebra with basis $\{1,e_{1},e_{2},e_{3}\}~$and $\mathbb{H(}1,1)$ be the
quaternion division algebra with basis $\{1,i,j,k\}$. The multiplication
table is given below:%
\begin{equation*}
\begin{tabular}{c||c|c|c|c|}
$\cdot $ & $1$ & $e_{1}$ & $e_{2}$ & $e_{3}$ \\ \hline\hline
$1$ & $1$ & $e_{1}$ & $e_{2}$ & $e_{3}$ \\ \hline
$e_{1}$ & $e_{1}$ & $-\gamma _{1}$ & $e_{3}$ & $-\gamma _{1}e_{2}$ \\ \hline
$e_{2}$ & $e_{2}$ & $-e_{3}$ & $-\gamma _{2}$ & $\gamma _{2}e_{1}$ \\ \hline
$e_{3}$ & $e_{3}$ & $\gamma _{1}e_{2}$ & $-\gamma _{2}e_{1}$ & $-\gamma
_{1}\gamma _{2}$ \\ \hline
\end{tabular}%
\end{equation*}

The algebra $\mathbb{H}(1,-1)$ is called \textit{the} \textit{algebra of
coquaternions} \cite{Cockle} or also called \textit{the algebra of
para-quaternions} \cite{Ivanov}, or \textit{the algebra of split-quaternions}
\cite{Erdogdu}, \cite{Ozdemir}, or \textit{the algebra of anti-quaternions},
or \textit{the algebra of pseudo-quaternions} \cite[p. 389]{Rosenfeld}, or
\textit{hyperbolic quaternions} \cite{Carmody}. We denote by $%
\{1,i_{1},i_{2},i_{3}\}$ the basis of coquaternion algebra. \medskip

\textbf{Proposition 1.} (\cite{Lam}, Proposition 1.1) \textit{The quaternion
algebra} $\mathbb{H}\left( \beta _{1},\beta _{2}\right) $ \textit{is
isomorphic with quaternion algebra} $\mathbb{H}\left( x^{2}\beta
_{1},y^{2}\beta _{2}\right) ,$ \textit{where} $x,y\in K^{\ast }.$\textit{\ }$%
\Box \medskip $

The real octonion division algebras are a non-associative and
non-commutative extension of the algebra of quaternions. Among all the real
division algebras, octonion algebra forms the largest normed division
algebra.\qquad\ \ \ \ \qquad\ \ \ \qquad \qquad \qquad \qquad

Let $\mathbb{O}(\alpha ,\beta ,\gamma )$ be a generalized octonion algebra
over $\mathbb{R},$ with basis $\{1,f_{1},...,f_{7}\}\medskip \,\,$and the
multiplication given in the following table:\medskip \medskip

{\footnotesize $%
\begin{tabular}{c||c|c|c|c|c|c|c|c|}
$\cdot $ & $1$ & $\,\,\,f_{1}$ & $\,\,\,\,\,f_{2}$ & $\,\,\,\,f_{3}$ & $%
\,\,\,\,f_{4}$ & $\,\,\,\,\,\,f_{5}$ & $\,\,\,\,\,\,f_{6}$ & $%
\,\,\,\,\,\,\,f_{7}$ \\ \hline\hline
$\,1$ & $1$ & $\,\,\,f_{1}$ & $\,\,\,\,f_{2}$ & $\,\,\,\,f_{3}$ & $%
\,\,\,\,f_{4}$ & $\,\,\,\,\,\,f_{5}$ & $\,\,\,\,\,f_{6}$ & $%
\,\,\,\,\,\,\,f_{7}$ \\ \hline
$\,f_{1}$ & $\,\,f_{1}$ & $-\alpha $ & $\,\,\,\,f_{3}$ & $-\alpha f_{2}$ & $%
\,\,\,\,f_{5}$ & $-\alpha f_{4}$ & $-\,\,f_{7}$ & $\,\,\,\alpha f_{6}$ \\
\hline
$\,f_{2}$ & $\,f_{2}$ & $-f_{3}$ & $-\,\beta $ & $\,\,\beta f_{1}$ & $%
\,\,\,\,f_{6}$ & $\,\,\,\,\,f_{7}$ & $-\beta f_{4}$ & $-\beta f_{5}$ \\
\hline
$f_{3}$ & $f_{3}$ & $\alpha f_{2}$ & $-\beta f_{1}$ & $-\alpha \beta $ & $%
\,\,\,\,f_{7}$ & $-\alpha f_{6}$ & $\,\,\,\beta f_{5}$ & $-\alpha \beta
f_{4} $ \\ \hline
$f_{4}$ & $f_{4}$ & $-f_{5}$ & $-\,f_{6}$ & $-\,\,f_{7}$ & $-\,\gamma $ & $%
\,\,\,\gamma f_{1}$ & $\,\,\gamma f_{2}$ & $\,\,\,\,\,\gamma f_{3}$ \\ \hline
$\,f_{5}$ & $\,f_{5}$ & $\alpha f_{4}$ & $-\,f_{7}$ & $\,\alpha f_{6}$ & $%
-\gamma f_{1}$ & $-\,\alpha \gamma $ & $-\gamma f_{3}$ & $\,\alpha \gamma
f_{2}$ \\ \hline
$\,\,f_{6}$ & $\,\,f_{6}$ & $\,\,\,\,f_{7}$ & $\,\,\beta f_{4}$ & $-\,\beta
f_{5}$ & $-\gamma f_{2}$ & $\,\,\,\gamma f_{3}$ & $-\beta \gamma $ & $-\beta
\gamma f_{1}$ \\ \hline
$\,\,f_{7}$ & $\,\,f_{7}$ & $-\alpha f_{6}$ & $\,\beta f_{5}$ & $\alpha
\beta f_{4}$ & $-\gamma f_{3}$ & $-\alpha \gamma f_{2}$ & $\beta \gamma
f_{1} $ & $-\alpha \beta \gamma $ \\ \hline
\end{tabular}
\medskip $ }\medskip

The algebra $\mathbb{O}(\alpha ,\beta ,\gamma )$ is non-commutative,
non-associative but it is \textit{alternative}\ (i.e.\thinspace \thinspace $%
x^{2}y=x\left( xy\right) $ and $yx^{2}=\left( yx\right) x,\forall x,y\in
\mathbb{O}(\alpha ,\beta ,\gamma )$),\thinspace \textit{\ flexible}\ (i.e. $%
x\left( yx\right) =\left( xy\right) x,\forall x,y\in \mathbb{O}(\alpha
,\beta ,\gamma )$), \textit{power-associative}\ (i.e. for each $x\in \mathbb{%
O}(\alpha ,\beta ,\gamma )$ the subalgebra generated by $x$ is an
associative algebra).

If $a\in \mathbb{O}(\alpha ,\beta ,\gamma ),$ $%
a=a_{0}+a_{1}f_{1}+a_{2}f_{2}+a_{3}f_{3}+a_{4}f_{4}+a_{5}f_{5}+a_{6}f_{6}+a_{7}f_{7}
$ then $\bar{a}%
=a_{0}-a_{1}f_{1}-a_{2}f_{2}-a_{3}f_{3}-a_{4}f_{4}-a_{5}f_{5}-a_{6}f_{6}-a_{7}f_{7}
$ is called the \textit{conjugate} of the element $a.$ Let $A=\mathbb{O}%
(\alpha ,\beta ,\gamma )$ and $a\in A.$ We have that $t\left( a\right) =a+%
\overline{a}\in K$ and
\begin{equation*}
\,N\left( a\right) =a\overline{a}=a_{0}^{2}+\alpha a_{1}^{2}+\beta
a_{2}^{2}+\alpha \beta a_{3}^{2}+\gamma a_{4}^{2}+\alpha \gamma
a_{5}^{2}+\beta \gamma a_{6}^{2}+\alpha \beta \gamma a_{7}^{2}\in K.
\end{equation*}
These elements are called the \textit{trace}, respectively, the \textit{norm}
of the element $a\in $ $A.$ \thinspace \thinspace \thinspace It\thinspace
\thinspace \thinspace follows\thinspace \thinspace \thinspace that$%
\,\,\left( a+\overline{a}\right) a\,=a^{2}+\overline{a}a$=$a^{2}+n\left(
a\right) \cdot 1$ and\thinspace \thinspace $a^{2}-t\left( a\right) a+N\left(
a\right) =0,\forall a\in A,$ $\,$therefore the generalized octonion algebra
is \textit{quadratic}.

The subset $\,\,A_{0}$=$\{x\in A\,\,\,\mid \,t\left( x\right) $=$0\}$ of $%
A\,\,\,$is a subspace of the algebra $A$. It is obvious that $A=K\cdot
1\oplus A_{0},$ therefore each element $\,\,x\in \,\,A\,\,$\thinspace has
the form $x=x_{0}\cdot 1+\overrightarrow{x},\,$with $\,x_{0}\in K\,\,\,\,$and%
$\,\,\,\overrightarrow{x}\in A_{0}.$ For $K=\mathbb{R},$ we call $x_{0}$
\textit{the scalar part} and $\overrightarrow{x}$ \textit{the vector \ part}
for the octonion $x.$

If \thinspace for $x\in A,$ the relation $N\left( x\right) =0$ implies $x=0$%
, then the algebra $A$ is called a \textit{division} algebra. For other
details, the reader is referred to \cite{Schafer}.\bigskip

In the papers \cite{Flaut}, \cite{Flaut-St} are considered some algebraic
equations in generalized quaternion and octonion algebras. Due to
Proposition 1, in the present paper we reduced the study of an algebraic
equation in an arbitrary algebra $\mathbb{H}\left( \gamma _{1},\gamma
_{2}\right) $ with $\gamma _{1},\gamma _{2}\in \mathbb{R}\setminus \{0\}$ to
study of the corresponding algebraic equation in one of the following two
algebras: division quaternion algebra or coquaternion algebra. Moreover, De
Moivre's formula and Euler's formula\ in generalized quaternion algebras,
founded in \cite{Mamagani}, was proved using this new method, for $\gamma
_{1},\gamma _{2}>0.$ With this technique, the above mentioned results were
also obtained for the octonions.
\begin{equation*}
\end{equation*}

\textbf{1. An isomorphism between the algebras $\mathbb{H}(\gamma
_{1},\gamma _{2}),$with $\gamma _{1},\gamma _{2}>0,$ and $\mathbb{H}(1,1)$}%
\vspace{3mm}%
\begin{equation*}
\end{equation*}

Everywhere in this section, we will consider $\gamma _{1},\gamma _{2}>0$. An
isomorphism between the algebras $\mathbb{H}(\gamma _{1},\gamma _{2})$ and $%
\mathbb{H}(1,1)$ is given by the operator $A$ and its inverse $A^{-1}$, where%
\newline
\begin{equation*}
A:\quad e_{1}\mapsto i\sqrt{\gamma _{1}},\quad e_{2}\mapsto j\sqrt{\gamma
_{2}},\quad e_{3}\mapsto k\sqrt{\gamma _{1}\gamma _{2}}.
\end{equation*}

It is easy to prove the following properties for the operator $A:$ \vskip2mm
1) $A\left( \lambda x\right) =\lambda A\left( x\right) ,\forall \;\lambda
\in \mathbb{R},\forall \;x\in \mathbb{H}\left( \gamma _{1},\gamma
_{2}\right) ;$ \vskip2mm 2) $A\left( x+y\right) =A\left( x\right) +A\left(
y\right) ,\forall \;x,y\in \mathbb{H}\left( \gamma _{1},\gamma _{2}\right) ;$
\vskip2mm 3) $A\left( xy\right) =A\left( x\right) A\left( y\right) ,\forall
\;x,y\in \mathbb{H}\left( \gamma _{1},\gamma _{2}\right) .$ \vskip2mm From
here, it results that the operators $A$ and $A^{-1}$ are additive and
multiplicative. \vskip2mm

\textbf{Proposition 1.1.} \textit{The operators $A$ and $A^{-1}$} \textit{%
are continuous and their norms are equal with} $1$.\medskip

\textbf{Proof.} We denote by $\Vert \cdot \Vert _{\mathbb{H}(\gamma
_{1},\gamma _{2})}$ the Euclidian norm in $\mathbb{H}\left( \gamma
_{1},\gamma _{2}\right) $. Since the spaces $\mathbb{H}\left( \gamma
_{1},\gamma _{2}\right) $ and $\mathbb{H}(1,1)$ are normed spaces, then the
continuity of $A$ is equivalent with the boundedness of $A,$ i.e. there is a
real constant $c$ such that for all $x\in \mathbb{H}\left( \gamma
_{1},\gamma _{2}\right) ,$ we have
\begin{equation*}
\frac{\left\Vert A\left( x\right) \right\Vert _{\mathbb{H}\left( 1,1\right) }%
}{\left\Vert x\right\Vert _{\mathbb{H}\left( \gamma _{1},\gamma _{2}\right) }%
}\leq c.
\end{equation*}

It results that\newline
\begin{equation*}
\frac{\left\Vert x_{0}+x_{1}i\sqrt{\gamma _{1}}+x_{2}j\sqrt{\gamma _{2}}%
+x_{3}k\sqrt{\gamma _{1}\gamma _{2}}\right\Vert }{\left\Vert
x_{0}+x_{1}e_{1}+x_{2}e_{2}+x_{3}e_{3}\right\Vert }=
\end{equation*}%
\begin{equation*}
=\frac{\sqrt{x_{0}^{2}+x_{1}^{2}\gamma _{1}+x_{2}^{2}\gamma
_{2}+x_{3}^{2}\gamma _{3}}}{\sqrt{x_{0}^{2}+x_{1}^{2}\gamma
_{1}+x_{2}^{2}\gamma _{2}+x_{3}^{2}\gamma _{3}}}=1.
\end{equation*}%
$\,\Box \bigskip $

\textbf{2. The algebra $\mathbb{H}(\gamma _{1},\gamma _{2}),$ with }$\gamma
_{1},\gamma _{2}<0$ \textbf{\ or }$\gamma _{1}\gamma _{2}<0$\vspace{3mm}%
\begin{equation*}
\end{equation*}

In this situation, the algebra $\mathbb{H}(\gamma _{1},\gamma _{2})$ is
isomorphic with $\mathbb{H}(1,-1).$ We suppose first that $\gamma
_{1},\gamma _{2}<0$. An isomorphism between the algebra $\mathbb{H}(\gamma
_{1},\gamma _{2}),$ where $\gamma _{1},\gamma _{2}<0,$ and the algebra $%
\mathbb{H}(1,-1)$ is given by the operator $B$ and its inverse $B^{-1}$,
where\newline
\begin{equation*}
B:\quad e_{1}\mapsto i_{3}\sqrt{-\gamma _{1}},\quad e_{2}\mapsto i_{2}\sqrt{%
-\gamma _{2}},\quad e_{3}\mapsto i_{1}\sqrt{\gamma _{1}\gamma _{2}}.
\end{equation*}

For $\gamma _{1}>0,\gamma _{2}<0$, an isomorphism between the algebra $%
\mathbb{H}(\gamma _{1},\gamma _{2})$ and the algebra $\mathbb{H}(1,-1)$ is
given by the operator $C$ and its inverse $C^{-1}$, where\newline
\begin{equation*}
C:\quad e_{1}\mapsto i_{1}\sqrt{\gamma _{1}},\quad e_{2}\mapsto i_{2}\sqrt{%
-\gamma _{2}},\quad e_{3}\mapsto i_{3}\sqrt{-\gamma _{1}\gamma _{2}}.
\end{equation*}

For $\gamma _{1}<0,\gamma _{2}>0,$ an isomorphism between the algebra $%
\mathbb{H}(\gamma _{1},\gamma _{2})$ and the algebra $\mathbb{H}(1,-1)$ is
given by the operator $D$ and its inverse $D^{-1}$, where\newline
\begin{equation*}
D:\quad e_{1}\mapsto i_{3}\sqrt{-\gamma _{1}},\quad e_{2}\mapsto i_{1}\sqrt{%
\gamma _{2}},\quad e_{3}\mapsto i_{2}\sqrt{-\gamma _{1}\gamma _{2}}.
\end{equation*}

The properties of the operators $B,B^{-1},C,C^{-1},D,D^{-1}$ are similarly
with the properties of the operator $A$.

Since each algebra $\mathbb{H}(\gamma _{1},\gamma _{2})$~ is isomorphic with
algebra of quaternions or coquaternions, it results that the above operators
provide us a simple way to generalize known results in these two algebras to
generalized quaternion algebra.

\vspace{3mm}%
\begin{equation*}
\end{equation*}

\textbf{3. Application to the algebraic equations}\vspace{3mm}%
\begin{equation*}
\end{equation*}

Let $x=x_{0}+x_{1}e_{1}+x_{2}e_{3}+x_{3}e_{3}\in \mathbb{H}(\gamma
_{1},\gamma _{2})$ and let $f:\mathbb{H}(\gamma _{1},\gamma _{2})\rightarrow
\mathbb{H}(\gamma _{1},\gamma _{2})$ be a continuous function of the form $%
f(x)=f_{0}(x_{0},x_{1},x_{2},x_{3})+f_{1}(x_{0},x_{1},x_{2},x_{3})e_{1}+f_{2}(x_{0},x_{1},x_{2},x_{3})e_{2}+f_{3}(x_{0},x_{1},x_{2},x_{3})e_{3}
$. Let $F$ be the one of the operators $A,B,C$ or $D$, depending on the
signs of $\gamma _{1}$ and $\gamma _{2}$. We define the operator $\mathfrak{F%
}$ which for any continuous function $f,$ taking values in $\mathbb{H}%
(\gamma _{1},\gamma _{2}),$ maps it in the continuous function $\mathfrak{F}%
f,$ taking values in $\mathbb{H}(1,1)$ or $\mathbb{H}(1,-1)$ by the rule:%
\newline
\begin{equation*}
\mathfrak{F}f:=f_{0}+f_{1}\,F(e_{1})+f_{2}\,F(e_{2})+f_{3}\,F(e_{3}).%
\medskip
\end{equation*}

\textbf{Theorem 3.1.} \textit{Let }$x^{0}\in \mathbb{H}(\gamma _{1},\gamma
_{2})$\textit{\ be a root of the equation }$f(x)=0$\textit{$~$ in }$\mathbb{H%
}(\gamma _{1},\gamma _{2})$\textit{. Then }$F(x^{0})$\textit{\ is a root of
the equation $\mathfrak{F}$}$f\big(F(x)\big)=0$\textit{\ in }$\mathbb{H}(1,1)
$\textit{\ or }$\mathbb{H}(1,-1)$\textit{, depending on the signs of }$%
\gamma _{1}$\textit{\ and }$\gamma _{2}$\textit{. The converse is also true.}
\medskip

\textbf{Proof.} Let $\gamma _{1},\gamma _{2}>0$. Applying operator $A$ to
the equality $f(x^{0})=0$ and using the continuity of $A$, we obtain\newline
\begin{equation*}
A\big(f(x^{0})\big )=Af\big(A(x^{0})\big)=A(0)=0.\
\end{equation*}

To prove the converse statement we apply the operator $A^{-1}$ to the
equality $f(x^{0})=0$. The remaining cases can be proved
similarly.\thinspace $\Box $ \medskip

Therefore, all results from quaternionic equations and from coquaternionic
equations can be generalized in $\mathbb{H}(\gamma _{1},\gamma _{2}).$

It is known that any polynomial of degree $n$ with coefficients in a field \
$K$ has at most $n$ roots in $K$. If the coefficients are in $\mathbb{H}%
(1,1),$ the situation is different. For the real division quaternion algebra
over the real field, there is a kind of a fundamental theorem of algebra:
\textit{If a polynomial has only one term of the greatest degree, then it}
\textit{has at least one root in} $\mathbb{H}(1,1)$.  (\cite{Smith}, Theorem
65; \cite{Eilenberg}, Theorem 1).

We consider the polynomial of degree $n$ of the form
\begin{equation}
f(x)=a_{0}xa_{1}x\ldots a_{n-1}xa_{n}+\varphi (x),  \tag{3.1.}
\end{equation}%
where $x,\,a_{0},\,a_{1},\ldots,a_{n-1},\,a_{n}\in \mathbb{H}(1,1),$ with $%
a_{\ell }\neq 0$ for $\ell \in \{0,1,\ldots,n\}$ and $\varphi (x)$ is a sum
of a finite number of monomials of the form $b_{0}xb_{1}x\ldots
b_{t-1}xb_{t} $ where $t<n$. From the above, it results that the equation $%
f(x)=0$ has at least one root. Applying operator $A^{-1}$ to this last
equality, the equation $\big(A^{-1}f\big)\big(A^{-1}(x)\big)=0,$ with $%
x=x_{0}+x_{1}e_{1}+x_{2}e_{2}+x_{3}e_{3},$ has at least one root in $\mathbb{%
H}(\gamma _{1},\gamma _{2})$. Therefore, we proved the following
result:\medskip

\textbf{Theorem 3.2.} \textit{In the generalized quaternion algebra }$%
\mathbb{H}(\gamma _{1},\gamma _{2})$\textit{\ with }$\gamma _{1},\gamma
_{2}>0$\textit{\ any polynomial of the form }$(3.1)$\textit{\ has at least
one root.}$\,\Box $\medskip

In the following, we consider the equation $x^{2}+ax+b=0,$ where $\mathbb{H}%
(\gamma _{1},\gamma _{2})$ with $\gamma _{1},\gamma _{2}>0$. We say that a
root $x_{0}=a_{0}+a_{1}i_{1}+a_{2}i_{2}+a_{3}i_{3}\in \mathbb{H}(\gamma
_{1},\gamma _{2})$,\ $\gamma _{1},\gamma _{2}>0$ is a \textit{ellipsoidal
root} if every element $x_{1}\in \mathbb{H}(\gamma _{1},\gamma _{2})$ of the
form $x_{1}=a_{0}+b_{1}i_{1}+b_{2}i_{2}+b_{3}i_{3}$ such that $\gamma
_{1}a_{1}^{2}+\gamma _{2}a_{2}^{2}+\gamma _{1}\gamma _{2}a_{3}^{2}=\gamma
_{1}b_{1}^{2}+\gamma _{2}b_{2}^{2}+\gamma _{1}\gamma _{2}b_{3}^{2}$ is also
a root of this equation.

Using Theorem 3 from \cite{Mierz} and Theorem 3.1, we just proved the
following theorem:\medskip

\textbf{Theorem 3.3.} \textit{In the generalized quaternion algebra }$%
\mathbb{H}(\gamma _{1},\gamma _{2})$\textit{\ with }$\gamma _{1},\gamma
_{2}>0$\textit{\ the equation }$x^{2}+ax+b=0$\textit{\ has ellipsoidal root
if and only if }$a$\textit{\ and }$b$\textit{\ are real numbers and }$%
a^{2}-4b<0$\textit{.}\ $\Box $\medskip

More results on the structure of roots of the quadratic quaternionic
equations can be found in see \cite{Szpakowski}, \cite{Mierz-1}, \cite%
{Mierz-2}, \cite{Niven}.

In the following, we apply the above results to the coquaternion algebra. We
will consider one of the three cases: $\gamma _{1},\gamma _{2}<0$;\thinspace
$\gamma _{1}<0,\gamma _{2}>0$;\ $\gamma _{1}>0,\gamma _{2}<0$. We say that a
root $x_{0}=a_{0}+a_{1}i_{1}+a_{2}i_{2}+a_{3}i_{3}\in \mathbb{H}(\gamma
_{1},\gamma _{2})$ is a \textit{hyperboloidal root} if every element $%
x_{1}\in \mathbb{H}(\gamma _{1},\gamma _{2})$ of the form $%
x_{1}=a_{0}+b_{1}i_{1}+b_{2}i_{2}+b_{3}i_{3}$ such that $\gamma
_{1}a_{1}^{2}+\gamma _{2}a_{2}^{2}+\gamma _{1}\gamma _{2}a_{3}^{2}=\gamma
_{1}b_{1}^{2}+\gamma _{2}b_{2}^{2}+\gamma _{1}\gamma _{2}b_{3}^{2}$ is also
a root of this equation. Using Theorem 2.5 of \cite{Pogoruy} and the above
Theorem 3.1, we proved:\medskip

\textbf{Theorem 3.4.} \textit{Suppose that coefficients }$r_{k}$, $%
k=0,1,\ldots ,n$\textit{\ of the polynomial equation }$%
r_{n}x^{n}+r_{n-1}x^{n-1}+\ldots +r_{0}=0$, \ $x\in \mathbb{H}(\gamma
_{1},\gamma _{2})$\textit{\ are real numbers. Therefore all its non-real
roots are hyperboloidal.}\thinspace $\Box $ \medskip

Solutions of linear equations in coquaternionic algebra can be found in \cite%
{Janovska}, \cite{Erdogdu} and solutions of linear equations in quaternion
algebra can be found, for example, in \cite{Tian}, \cite{Shpakivskyi}.
\medskip

Let $x=x_{0}+x_{1}i+x_{2}j+x_{3}k\in \mathbb{H}(1,1)$, \ $%
x_{0},x_{1},x_{2},x_{3}\in \mathbb{R}.$ We denoted by\newline
\begin{equation*}
\widetilde{x}=A^{-1}(x)=x_{0}+x_{1}\frac{e_{1}}{\sqrt{\gamma _{1}}}%
+x_{2}e_{2}\frac{e_{2}}{\sqrt{\gamma _{2}}}+x_{3}e_{3}\frac{e_{1}}{\sqrt{%
\gamma _{1}\gamma _{2}}}.\medskip
\end{equation*}

\textbf{Example 3.5.} Consider the general polynomial equation in $\mathbb{H}%
(1,1):$\newline
\begin{equation*}
\underset{p=1}{\overset{n}{\sum }}\left( \underset{l=1}{\overset{m_{p}}{\sum
}}a_{p,l,1}xa_{p,l,2}x...a_{p,l,p}xa_{p,l,p+1}\right) +c=0.
\end{equation*}

This equation is "equivalent" to the following equation in $\mathbb{H}%
(\gamma _{1},\gamma _{2}),$ with $\gamma _{1},\gamma _{2}>0$:\newline
\begin{equation*}
\underset{p=1}{\overset{n}{\sum }}\left( \underset{l=1}{\overset{m_{p}}{\sum
}}\widetilde{a}_{p,l,1}\widetilde{x}\,\widetilde{a}_{p,l,2}\,\widetilde{x}%
\ldots \widetilde{a}_{p,l,p}\,\widetilde{x}\,\widetilde{a}_{p,l,p+1}\right) +%
\widetilde{c}=0.
\end{equation*}

Therefore, the algebraic equations $f(x)=0,$ with $f$ continuous, in usual
quaternions $\mathbb{H}(1,1)$ can be reduced to the similar equations in all
algebras $\mathbb{H}(\gamma _{1},\gamma _{2})$ for $\gamma _{1},\gamma
_{2}>0 $ and vice versa. A similar result can be found for different signs
of $\gamma _{1}$ and $\gamma _{2}$.

\begin{equation*}
\end{equation*}

\textbf{4. De Moivre's formula}\vspace{3mm}%
\begin{equation*}
\end{equation*}

In the following, we will use some ideas and notations from \cite{Cho}. Let $%
q=q_{0}+q_{1}e_{1}+q_{2}e_{2}+q_{3}e_{3}\in \mathbb{H}(\gamma _{1},\gamma
_{2})$,\ $\gamma _{1},\gamma _{2}>0$, \ $q_{0},q_{1},q_{2},q_{3}\in \mathbb{R%
}$ and
\begin{equation*}
|q|=\sqrt{q_{0}^{2}+\gamma _{1}q_{1}^{2}+\gamma _{2}q_{2}^{2}+\gamma
_{1}\gamma _{2}q_{3}^{2}}.
\end{equation*}
Consider the sets\newline
\begin{equation*}
\mathcal{S}_{G}^{3}=\{q\in \mathbb{H}(\gamma _{1},\gamma _{2}),\gamma
_{1},\gamma _{2}>0:\,|q|=1\},\newline
\end{equation*}%
\begin{equation*}
\mathcal{S}_{G}^{2}=\{q\in \mathbb{H}\left( \gamma _{1},\gamma _{2}\right)
,\gamma _{1},\gamma _{2}>0:\,q_{0}=0,|q|=1\}.
\end{equation*}

Any $q\in \mathcal{S}_{G}^{3}$ can be expressed as $q=\cos \theta
+\varepsilon \sin \theta ,$ where\newline
\begin{equation*}
\cos \theta =q_{0},\;\;\varepsilon =\frac{q_{1}e_{1}+q_{2}e_{2}+q_{3}e_{3}}{%
\sqrt{\gamma _{1}q_{1}^{2}+\gamma _{2}q_{2}^{2}+\gamma _{1}\gamma
_{2}q_{3}^{2}}}.
\end{equation*}

Using Proposition 2 from \cite{Cho} and applying the operator $A^{-1}$ we
will find De Moivre's formula for $\mathbb{H}\left( \gamma _{1},\gamma
_{2}\right) ,\gamma _{1},\gamma _{2}>0$.\medskip

\textbf{Theorem 4.1.} \textit{Let }$q=\cos \theta +\varepsilon \sin \theta
\in $\textit{$\mathcal{S}$}$_{G}^{3}$,\ $\theta \in \mathbb{R}$\textit{.
Then }$q^{n}=\cos n\theta +\varepsilon \sin n\theta $\textit{\ for every
integer }$n$\textit{.}\medskip

Theorem 4.1 is the same with Theorem 7 from the paper \cite{Mamagani},
obtained with another proof.$\Box \medskip $

Using Corollary 3 from \cite{Cho} and Theorem 3.1, we obtain the next
statement.\medskip

\textbf{Proposition 4.2.} \textit{\ i) In }$\mathbb{H}(\gamma _{1},\gamma
_{2}),$ $\gamma _{1},\gamma _{2}>0$\textit{\ the equation }$x^{n}=1$\textit{%
\ whit }$n$\textit{\ integer and }$n\geq 3$\textit{\ has infinity of roots,
namely }
\begin{equation*}
q=\cos \frac{2\pi }{n}+\varepsilon \sin \frac{2\pi }{n}\in \mathcal{S}%
_{G}^{3},\ \varepsilon \in \mathcal{S}_{G}^{2}.
\end{equation*}

\textit{\ ii) In }$\mathbb{H}(\gamma _{1},\gamma _{2}),$ $\gamma _{1},\gamma
_{2}>0$\textit{\ the equation }$x^{n}=a$, $n\in \mathbb{N}$, $a\in \mathbb{R}
$\textit{\ has infinity of roots, namely }$\sqrt[n]{a}\,q,$\textit{\ where }$%
q=\cos \frac{2\pi }{n}+\varepsilon \sin \frac{2\pi }{n}\in \mathcal{S}%
_{G}^{3},$\textit{\ \thinspace with\thinspace\ }$\varepsilon \in $\textit{$%
\mathcal{S}_{G}^{2}.$ If }$n$\textit{\ is even it is necessary that }$a>0$%
\textit{.\medskip }$\Box \medskip $%
\begin{equation*}
\end{equation*}

\textbf{5.} \bigskip \textbf{De Moivre's formula and Euler's formula for
octonions}

\begin{equation*}
\end{equation*}

\

In the following, we will generalize in a natural way De Moivre formula and
Euler's formula for the division octonion algebra \ $\mathbb{O}\left(
1,1,1\right) .$ For this, we will use some ideas and notations from \cite%
{Cho}. We consider the sets
\begin{equation*}
\mathcal{S}^{3}=\{a\in \mathbb{O}(1,1,1):\;N(a)=1\},
\end{equation*}%
\begin{equation*}
\mathcal{S}_{G}^{3}=\{a\in \mathbb{O}(\alpha ,\beta ,\gamma ):\;N(a)=1\},
\end{equation*}%
\begin{equation*}
\mathcal{S}^{2}=\{a\in \mathbb{O}(1,1,1):\;t(a)=0,N\left( a\right) =1\}.
\end{equation*}

\begin{equation*}
\mathcal{S}_{G}^{2}=\{a\in \mathbb{O}( 1,1,1):\;t(a) =0,N\left( a\right)
=1\}.
\end{equation*}%
We remark that for all elements $a\in \mathcal{S}^{2},$ we have $a^{2}=-1.$
Let $a\in \mathcal{S}%
^{3},a=a_{0}+a_{1}f_{1}+a_{2}f_{2}+a_{3}f_{3}+a_{4}f_{4}+a_{5}f_{5}+a_{6}f_{6}+a_{7}f_{7}.
$ This element can be write under the form%
\begin{equation*}
a=\cos \lambda +w\sin \lambda ,
\end{equation*}%
where $\cos \lambda =a_{0}~$ and
\begin{equation*}
w=\frac{
a_{1}f_{1}+a_{2}f_{2}+a_{3}f_{3}+a_{4}f_{4}+a_{5}f_{5}+a_{6}f_{6}+a_{7}f_{7}
}{\sqrt{%
a_{1}^{2}+a_{2}^{2}+a_{3}^{2}+a_{4}^{2}+a_{5}^{2}+a_{6}^{2}+a_{7}^{2} }}=
\end{equation*}
\begin{equation*}
=\frac{
a_{1}f_{1}+a_{2}f_{2}+a_{3}f_{3}+a_{4}f_{4}+a_{5}f_{5}+a_{6}f_{6}+a_{7}f_{7}
}{\sqrt{1-a_{0}^{2}}}.
\end{equation*}

Since $w^{2}=-1,$ we obtain the following Euler's formula:%
\begin{eqnarray*}
e^{\lambda w} &=&\overset{\infty }{\underset{i=1}{\sum }}\frac{\left(
\lambda w\right) ^{n}}{n!}=\overset{\infty }{\underset{i=1}{\sum }}\frac{%
\left( -1\right) ^{n}\lambda ^{2n}}{\left( 2n\right) !}+w\overset{\infty }{%
\underset{i=1}{\sum }}\frac{\left( -1\right) ^{n-1}\lambda ^{2n-1}}{\left(
2n-1\right) !}= \\
&=&\cos \lambda +w\sin \lambda .
\end{eqnarray*}

\textbf{Proposition 5.1.} \textit{The cosinus function is constant for all
elements in} $\mathcal{S}^{2}.\medskip $

\textbf{Proof.} Indeed, $\cos w=\overset{\infty }{\underset{i=1}{\sum }}%
\frac{\left( -1\right) ^{n}\lambda ^{2n}}{\left( 2n\right) !}=\cos i.\Box
\medskip $

\textbf{Proposition 5.2.} For $w\in \mathcal{S}^{2},$ we have $(\cos \lambda
_{1}+w\sin \lambda _{1})(\cos \lambda _{2}+w\sin \lambda _{2})=\cos (\lambda
_{1}+\lambda _{2})+w\sin (\lambda _{1}+\lambda _{2}).\medskip $\qquad

\textbf{Proof.} By straightforward calculations $\Box \medskip $

\textbf{Proposition 5.3.} (De Moivre formula for octonions) \ \textit{With
the above notations, we have that}
\begin{equation*}
a^{n}=e^{n\lambda w}=\left( \cos \lambda +w\sin \lambda \right) ^{n}=\cos
n\lambda +w\sin n\lambda ,
\end{equation*}%
\textit{where} $a\in \mathcal{S}^{3},$ $n\in \mathbb{Z}$ \textit{and} $%
\lambda \in \mathbb{R}.\medskip $

\textbf{Proof.} For $n>0,$ by induction. We obtain
\begin{eqnarray*}
a^{n+1} &=&\left( \cos \lambda +w\sin \lambda \right) ^{n+1}= \\
&=&\left( \cos \lambda +w\sin \lambda \right) ^{n}\left( \cos \lambda +w\sin
\lambda \right) = \\
&=&(\cos n\lambda +w\sin n\lambda )\left( \cos \lambda +w\sin \lambda
\right) = \\
&=&\cos (n+1)\lambda +w\sin (n+1)\lambda .
\end{eqnarray*}

Since $a^{-1}=\cos \lambda -w\sin \lambda =\cos (-\lambda )+w\sin (-\lambda
),$ it results the asked formula for all $n\in \mathbb{Z}.\Box \medskip $

\textbf{Remark 5.4.} It is known that any polynomial of degree $n$ with
coefficients in a field \ $K$ has at most $n$ roots in $K$. If the
coefficients are in $\mathbb{O}\left( 1,1,1\right) $ there is a kind of a
fundamental theorem of algebra: If a polynomial has only one term of the
greatest degree, then it has at least one root in $\mathbb{O}\left(
1,1,1\right) $ (see \cite{Smith}, Theorem 65). \medskip

\textbf{Theorem 5.5.} \textit{Equation} $x^{n}=a,~$\textit{where} $a\in
\mathbb{O}\left( 1,1,1\right) \setminus\mathbb{R},$ \textit{has }$n$ \textit{%
\ roots.}\medskip

\textbf{Proof.} The octonion $a$ can be written under the form $a=\sqrt{%
N\left( a\right) }\frac{a}{\sqrt{N\left( a\right) }}.$ The octonion $b=\frac{%
a}{\sqrt{N\left( a\right) }}$ is in $\mathcal{S}^{3},$ then \ we can find
the elements $w\in \mathcal{S}^{2}$ and $\lambda \in \mathbb{R}$ such that $%
b=\cos \lambda +w\sin \lambda .$ From Proposition 5.3, we have that the
solutions of the above equation are $x_{r}=\sqrt[n]{Q}\left( \cos \frac{%
\lambda +2r\pi }{n}+w\sin \frac{\lambda +2r\pi }{n}\right) ,$ where $Q=\sqrt{%
N\left( a\right) }$ and $\ r\in \{0,1,...,n-1\}.\Box \medskip $

\textbf{Corollary 5.6.} \textit{If} $a\in \mathbb{R},$ \textit{therefore the
equation} $x^{n}=a$ \textit{has an infinity of roots.\medskip }

\textbf{Proof.} Indeed, if $a\in \mathbb{R}$, we can write $a=a\cdot
1=a\left( \cos 2\pi +w\sin 2\pi \right) ,$ where $w\in \mathcal{S}^{2}$ is
an arbitrary element. $\Box \medskip $

\textbf{Remark 5.7.} The rotation of the octonion $x\in \mathbb{O}\left(
1,1,1\right) $ on the angle $\lambda $ around the unit vector $w\in \mathcal{%
S}^{2}$\ is defined by the formula
\begin{equation*}
x^{r}=\overline{u}xu,
\end{equation*}%
where $u\in \mathcal{S}^{3},u=\cos \frac{\lambda }{2}+w\sin \frac{\lambda }{2%
}$ and $\overline{u}=\cos \frac{\lambda }{2}-w\sin \frac{\lambda }{2}.$

Using the form $x=x_{0}\cdot 1+\overrightarrow{x},$ $y=y_{0}\cdot 1+%
\overrightarrow{y}\ \ \ $for the octonions $x,y\in \mathbb{O}\left(
1,1,1\right) ,$ we obtain the following expression for the product of two
octonions:%
\begin{eqnarray*}
xy &=&\left( x_{0}\cdot 1+\overrightarrow{x}\right) \left( y_{0}\cdot 1+%
\overrightarrow{y}\right) = \\
&=&x_{0}y_{0}\cdot 1+x_{0}\overrightarrow{y}+y_{0}\overrightarrow{x}+<%
\overrightarrow{x},\overrightarrow{y}>+\overrightarrow{x}\times
\overrightarrow{y},
\end{eqnarray*}%
\newline
where $<\overrightarrow{x},\overrightarrow{y}>$ is the inner product of two
octonionic-vector and $\overrightarrow{x}\times \overrightarrow{y}$ is the
cross product. From here, we obtain that $x^{r}=\overline{u}xu=$\newline
$=\left( \cos \frac{\lambda }{2}-w\sin \frac{\lambda }{2}\right) \left(
x_{0}\cdot 1+\overrightarrow{x}\right) \left( \cos \frac{\lambda }{2}+w\sin
\frac{\lambda }{2}\right) =$\newline
$=x_{0}+\overrightarrow{x}\cos \lambda -w<w,\overrightarrow{x}>\left( 1-\cos
\lambda \right) -\left( w\times \overrightarrow{x}\right) \sin \lambda .$ It
results that that rotation does not transform the octonion-scalar part, but
the octonion-vector part $\overrightarrow{x}\ $is rotated on the angle $%
\lambda $ around $w$. \newline

\begin{equation*}
\end{equation*}

\textbf{6. An isomorphism between the algebras $\mathbb{O}%
(\alpha,\beta,\gamma)$ and some its applications} \vspace{3mm}
\begin{equation*}
\end{equation*}

In the following, we will consider the generalized real octonion algebra $%
\mathbb{O}(\alpha ,\beta ,\gamma )$ and the algebras $\mathbb{O}( 1,1,1)$
and $\mathbb{O}( 1,1,-1)$. Let $\{1,f_{1},\ldots,f_{7}\}$ be a basis in $%
\mathbb{O} (\alpha ,\beta ,\gamma )$, and $\{1,\widetilde{f_{1}},\ldots,%
\widetilde{f_{7}}\}$ be the canonical basis in $\mathbb{O}( 1,1,1)$, and $%
\{1,\widehat{f_{1}},\ldots,\widehat{f_{7}}\}$ be the canonical basis in $%
\mathbb{O} (1,1,-1)$.

We prove that the algebra $\mathbb{O}(\alpha ,\beta ,\gamma )$ with $\alpha
,\beta ,\gamma \in \mathbb{R}\setminus \{0\}$ is isomorphic with algebra $%
\mathbb{O}(1,1,1)$ or $\mathbb{O}(1,1,-1)$ and indicate the formulae to pass
from one basis to another basis. Thus, if $\alpha ,\beta ,\gamma >0$ then
the real octonion algebra $\mathbb{O}(\alpha ,\beta ,\gamma )$ is isomorphic
with algebra $\mathbb{O}(1,1,1)$ and this isomorphism is given by the
relations:
\begin{eqnarray*}
A_{1}:\quad  &f_{1}\mapsto \widetilde{f_{1}}\sqrt{\alpha },&\;\;f_{2}\mapsto
\widetilde{f_{2}}\sqrt{\beta },\;\;f_{3}\mapsto \widetilde{f_{3}}\sqrt{%
\alpha \beta }, \\
&f_{4}\mapsto \widetilde{f_{4}}\sqrt{\gamma },&\;\;f_{5}\mapsto \widetilde{%
f_{5}}\sqrt{\alpha \gamma },\;\;f_{6}\mapsto \widetilde{f_{6}}\sqrt{\beta
\gamma },\;\;f_{7}\mapsto \widetilde{f_{7}}\sqrt{\alpha \beta \gamma }\,.
\end{eqnarray*}

If $\alpha ,\beta >0,\gamma <0,$ then the real octonion algebra $\mathbb{O}%
(\alpha ,\beta ,\gamma )$ is isomorphic with algebra $\mathbb{O}(1,1,-1)$
and this isomorphism is given by the relations:
\begin{eqnarray*}
A_{2}:\quad  &f_{1}\mapsto \widehat{f_{1}}\sqrt{\alpha },&\;\;f_{2}\mapsto
\widehat{f_{2}}\sqrt{\beta },\;\;f_{3}\mapsto \widehat{f_{3}}\sqrt{\alpha
\beta }, \\
&f_{4}\mapsto \widehat{f_{4}}\sqrt{-\gamma },&\;\;f_{5}\mapsto \widehat{f_{5}%
}\sqrt{-\alpha \gamma },\;\;f_{6}\mapsto \widehat{f_{6}}\sqrt{-\beta \gamma }%
,\;\;f_{7}\mapsto \widehat{f_{7}}\sqrt{-\alpha \beta \gamma }\,.
\end{eqnarray*}

If $\alpha ,\gamma >0,\beta <0,$ then the real octonion algebra $\mathbb{O}%
(\alpha ,\beta ,\gamma )$ is isomorphic with algebra $\mathbb{O}(1,1,-1)$
and this isomorphism is given by the relations:
\begin{eqnarray*}
A_{3}:\quad  &f_{1}\mapsto \widehat{f_{1}}\sqrt{\alpha },&\;\;f_{2}\mapsto
\widehat{f_{4}}\sqrt{-\beta },\;\;f_{3}\mapsto \widehat{f_{5}}\sqrt{-\alpha
\beta }, \\
&f_{4}\mapsto \widehat{f_{2}}\sqrt{\gamma },&\;\;f_{5}\mapsto \widehat{f_{3}}%
\sqrt{\alpha \gamma },\;\;f_{6}\mapsto \widehat{f_{6}}\sqrt{-\beta \gamma }%
,\;\;f_{7}\mapsto \widehat{f_{7}}\sqrt{-\alpha \beta \gamma }\,.
\end{eqnarray*}

If $\alpha >0,\beta ,\gamma <0,$ then the real octonion algebra $\mathbb{O}%
(\alpha ,\beta ,\gamma )$ is isomorphic with algebra $\mathbb{O}(1,1,-1)$
and this isomorphism is given by the relations:
\begin{eqnarray*}
A_{4}:\quad  &f_{1}\mapsto \widehat{f_{1}}\sqrt{\alpha },&\;\;f_{2}\mapsto
\widehat{f_{4}}\sqrt{-\beta },\;\;f_{3}\mapsto \widehat{f_{5}}\sqrt{-\alpha
\beta }, \\
&f_{4}\mapsto \widehat{f_{6}}\sqrt{-\gamma },&\;\;f_{5}\mapsto \widehat{f_{7}%
}\sqrt{-\alpha \gamma },\;\;f_{6}\mapsto \widehat{f_{2}}\sqrt{\beta \gamma }%
,\;\;f_{7}\mapsto \widehat{f_{3}}\sqrt{\alpha \beta \gamma }\,.
\end{eqnarray*}

If $\alpha <0,\beta ,\gamma >0,$ then the real octonion algebra $\mathbb{O}%
(\alpha ,\beta ,\gamma )$ is isomorphic with algebra $\mathbb{O}(1,1,-1)$
and this isomorphism is given by the relations:
\begin{eqnarray*}
A_{5}:\quad  &f_{1}\mapsto \widehat{f_{4}}\sqrt{-\alpha },&\;\;f_{2}\mapsto
\widehat{f_{1}}\sqrt{\beta },\;\;f_{3}\mapsto \widehat{f_{5}}\sqrt{-\alpha
\beta }, \\
&f_{4}\mapsto \widehat{f_{2}}\sqrt{\gamma },&\;\;f_{5}\mapsto \widehat{f_{6}}%
\sqrt{-\alpha \gamma },\;\;f_{6}\mapsto \widehat{f_{3}}\sqrt{\beta \gamma }%
,\;\;f_{7}\mapsto \widehat{f_{7}}\sqrt{-\alpha \beta \gamma }\,.
\end{eqnarray*}

If $\alpha ,\gamma <0,\beta >0,$ then real octonion algebra $\mathbb{O}%
(\alpha ,\beta ,\gamma )$ is isomorphic with algebra $\mathbb{O}(1,1,-1)$
and this isomorphism is given  by the relations:
\begin{eqnarray*}
A_{6}:\quad  &f_{1}\mapsto \widehat{f_{4}}\sqrt{-\alpha },&\;\;f_{2}\mapsto
\widehat{f_{1}}\sqrt{\beta },\;\;f_{3}\mapsto \widehat{f_{5}}\sqrt{-\alpha
\beta }, \\
&f_{4}\mapsto \widehat{f_{6}}\sqrt{-\gamma },&\;\;f_{5}\mapsto \widehat{f_{2}%
}\sqrt{\alpha \gamma },\;\;f_{6}\mapsto \widehat{f_{7}}\sqrt{-\beta \gamma }%
,\;\;f_{7}\mapsto \widehat{f_{3}}\sqrt{\alpha \beta \gamma }\,.
\end{eqnarray*}

If $\alpha ,\beta <0,\gamma >0,$ then the real octonion algebra $\mathbb{O}%
(\alpha ,\beta ,\gamma )$ is isomorphic with algebra $\mathbb{O}(1,1,-1)$
and this isomorphism is given by the relations:
\begin{eqnarray*}
A_{7}:\quad  &f_{1}\mapsto \widehat{f_{4}}\sqrt{-\alpha },&\;\;f_{2}\mapsto
\widehat{f_{5}}\sqrt{-\beta },\;\;f_{3}\mapsto \widehat{f_{1}}\sqrt{\alpha
\beta }, \\
&f_{4}\mapsto \widehat{f_{2}}\sqrt{\gamma },&\;\;f_{5}\mapsto \widehat{f_{6}}%
\sqrt{-\alpha \gamma },\;\;f_{6}\mapsto \widehat{f_{7}}\sqrt{-\beta \gamma }%
,\;\;f_{7}\mapsto \widehat{f_{3}}\sqrt{\alpha \beta \gamma }\,.
\end{eqnarray*}

If $\alpha ,\beta ,\gamma <0,$ then the real octonion algebra $\mathbb{O}%
(\alpha ,\beta ,\gamma )$ is isomorphic with algebra $\mathbb{O}(1,1,-1)$
and this isomorphism is given  by the relations:
\begin{eqnarray*}
A_{8}:\quad  &f_{1}\mapsto \widehat{f_{4}}\sqrt{-\alpha },&\;\;f_{2}\mapsto
\widehat{f_{5}}\sqrt{-\beta },\;\;f_{3}\mapsto \widehat{f_{1}}\sqrt{\alpha
\beta }, \\
&f_{4}\mapsto \widehat{f_{6}}\sqrt{-\gamma },&\;\;f_{5}\mapsto \widehat{f_{2}%
}\sqrt{\alpha \gamma },\;\;f_{6}\mapsto \widehat{f_{3}}\sqrt{\beta \gamma }%
,\;\;f_{7}\mapsto \widehat{f_{7}}\sqrt{-\alpha \beta \gamma }\,.
\end{eqnarray*}

It is easy to prove that the operators $A_{k}$, \thinspace $k=\overline{1,8}$
is additive and multiplicative. The following statement can be proved
completely analogous as Proposition 1.1. $\medskip $

\textbf{Proposition 6.1.} \textit{The operators }$A_{k}$,\thinspace\ $k=%
\overline{1,8}$\textit{\ are continuous and have norm } $1$.$\Box $\medskip

Let $x=x_{0}+\sum\limits_{k=1}^{7}x_{k}f_{k}\in \mathbb{O}(\alpha ,\beta
,\gamma )$ and let $g:\mathbb{O}(\alpha ,\beta ,\gamma )\rightarrow \mathbb{O%
}(\alpha ,\beta ,\gamma )$ be a continuous function of the form $%
g(x)=g_{0}(x_{0},\ldots ,x_{7})+\sum\limits_{k=1}^{7}g_{k}(x_{0},\ldots
,x_{7})f_{k}$. Let $L$ be  one of the operators $A_{k}$,\thinspace\ $k=%
\overline{1,8}$, depending on the signs of $\alpha ,\beta $ and $\gamma $.
We define the operator $\mathfrak{L}$\thinspace\ by the rule:\newline
\begin{equation*}
\mathfrak{L}g:=f_{0}+\sum\limits_{k=1}^{7}g_{k}\,L(f_{k}).
\end{equation*}

The operator $\mathfrak{L}$ for any continuous function $g,$ taking values
in $\mathbb{O}(\alpha ,\beta ,\gamma ),$ maps it in the continuous function $%
\mathfrak{L}g,$ taking values in $\mathbb{O}(1,1,1)$ or $\mathbb{O}(1,1,-1),$

The following statement can be analogously proved  as Theorem 3.1.\medskip

\textbf{Theorem 6.2.} \textit{Let }$x^{0}\in \mathbb{O}(\alpha ,\beta
,\gamma ),$\textit{\ be a root of the equation }$g(x)=0$\textit{\ in }$%
\mathbb{O}(\alpha ,\beta ,\gamma )$\textit{. Then }$L(x^{0})$\textit{\ is a
root of the equation $\mathfrak{L}$}$g\big(L(x)\big)=0$\textit{\ in }$%
\mathbb{O}(1,1,1)$\textit{\ or }$\mathbb{O}(1,1,-1)$\textit{, depending on
the signs of }$\alpha ,\beta ,$\textit{\ and }$\gamma $\textit{. The
converse is also true.} $\Box $ \medskip

Thus, the study of algebraic equations in an arbitrary algebra $\mathbb{O}%
(\alpha ,\beta ,\gamma )$ with $\alpha ,\beta ,\gamma \in \mathbb{R}%
\setminus \{0\}$ was reduced to study of the corresponding algebraic
equation in one of the following two algebras: division octonion algebra $%
\mathbb{O}(1,1,1)$ or algebra $\mathbb{O}(1,1,-1)$.

Using the above notations, we can prove the following theorem.\medskip

\textbf{Theorem 6.3.} \textit{Equation }$x^{n}=a,$\textit{\ where }$a\in
\mathbb{O}(\alpha ,\beta ,\gamma )\setminus \mathbb{R},\alpha ,\beta ,\gamma
>0,$\textit{\ has} $n$\textit{\ roots.\medskip }

\textbf{Proof.} The octonion $b=\frac{a}{\sqrt{N\left( a\right) }}$ is in $%
\mathcal{S}_{G}^{3},$ then there are $w\in \mathcal{S}_{G}^{2},%
\,w=A_{1}^{-1}\left( \widetilde{w}\right) ,\widetilde{w}\in \mathcal{S}^{2}$
and $\lambda \in \mathbb{R}$ such that $b=\cos \lambda +\widetilde{w}\sin
\lambda .$ From Proposition 5.3, we have that the solutions of the above
equation are $x_{r}=A_{1}^{-1}\left( \widetilde{x}_{r}\right) =\sqrt[2n]{%
N\left( a\right) }\left( \cos \frac{\lambda +2r\pi }{n}+\widetilde{w}\sin
\frac{\lambda +2r\pi }{n}\right) ,$ where $r\in \{0,1,\ldots ,n-1\}$ and $%
\widetilde{x}_{r}$ is a solution of the equation $\widetilde{x}^{n}=%
\widetilde{a}$ in $\mathbb{O}(1,1,1).\Box \medskip $

\textbf{Remark 6.4. } Using \ the operator $A_{1}^{-1},~$the rotation of the
octonion $x\in \mathbb{O}(\alpha ,\beta ,\gamma )$ on the angle $\lambda $
around the unit vector $w\in \mathcal{S}_{G}^{2}$\ is defined by the formula
\begin{equation*}
x^{r}=\overline{u}xu,
\end{equation*}%
where $u\in \mathcal{S}_{G}^{3},w\in \mathcal{S}_{G}^{2},~u=\cos \frac{%
\lambda }{2}+w\sin \frac{\lambda }{2}$ and $\overline{u}=\cos \frac{\lambda
}{2}-w\sin \frac{\lambda }{2}.$

By straightforward calculations, it results that  rotation does not
transform the octonion-scalar part, but the octonion-vector part $%
\overrightarrow{x}\ $is rotated on the angle $\lambda $ around $w$.\medskip

\textbf{Example 6.5.} \ \ 1) Let $a\in \mathcal{S}^{3},$\newline
$a=\frac{\sqrt{2}}{2}+\frac{1}{\sqrt{14}}\widetilde{f}_{1}+\frac{1}{\sqrt{14}%
}\widetilde{f}_{2}+\frac{1}{\sqrt{14}}\widetilde{f}_{3}+\frac{1}{\sqrt{14}}%
\widetilde{f}_{4}+\frac{1}{\sqrt{14}}\widetilde{f}_{5}+\frac{1}{\sqrt{14}}%
\widetilde{f}_{6}+\frac{1}{\sqrt{14}}\widetilde{f}_{7},$ we have $\cos
\lambda =\frac{\sqrt{2}}{2},\sin \lambda =\frac{\sqrt{2}}{2}.$ It results
that $a=\cos \frac{\pi }{4}+v\sin \frac{\pi }{4},$ where $v=\frac{1}{\sqrt{7}%
}(\widetilde{f}_{1}+\widetilde{f}_{2}+\widetilde{f}_{3}+\widetilde{f}_{4}+%
\widetilde{f}_{5}+\widetilde{f}_{6}+\widetilde{f}_{7}).$ The vector $a$
corresponds to the rotation of the space $\mathbb{R}^{8}$on the angle $\frac{%
\pi }{2}$ around the vector $v=\left( \frac{1}{\sqrt{7}},\frac{1}{\sqrt{7}}%
,\ldots ,\frac{1}{\sqrt{7}}\right) \in \mathbb{R}^{7}$ written in the
canonical basis.

2) In the algebra $\mathbb{O}(2,4,7),\,$\ for the above element $a\in
\mathcal{S}^{3},$ we have\newline
$b$=$A^{-1}\left( a\right) $=$\frac{\sqrt{2}}{2}$+$\frac{1}{2\sqrt{7}}f_{1}$+%
$\frac{1}{2\sqrt{14}}f_{2}$+$\frac{1}{7\sqrt{2}}f_{3}$+$\frac{1}{4\sqrt{7}}%
f_{4}$+$\frac{1}{14}f_{5}$+$\frac{1}{14\sqrt{2}}f_{6}$+$\frac{1}{28}f_{7}\in
\mathcal{S}_{G}^{3}$ and corresponds to the rotation of the space $\mathbb{R}%
^{8}$on the angle $\frac{\pi }{2}$ around the vector \newline
$v=\left( \frac{1}{2\sqrt{7}},\frac{1}{2\sqrt{14}},\frac{1}{7\sqrt{2}},\frac{%
1}{4\sqrt{7}},\frac{1}{14},\frac{1}{14\sqrt{2}},\frac{1}{28}\right) \in
\mathbb{R}^{7}$written in the basis $\{f_{1},\ldots ,f_{7}\}.$
\begin{equation*}
\end{equation*}

\bigskip \textbf{Case when }$\alpha =\beta =1,\gamma =-1$
\begin{equation*}
\end{equation*}

In this case, the octonion algebra $\mathbb{O}(1,1,-1)$ is not a division
algebra (is a split algebra). The norm of an octonion $a\in \mathbb{O}%
(1,1,-1),$ $%
a=a_{0}+a_{1}f_{1}+a_{2}f_{2}+a_{3}f_{3}+a_{4}f_{4}+a_{5}f_{5}+a_{6}f_{6}+a_{7}f_{7},
$ in this situation, can be positive, zero or negative. In the following, we
used \ definitions and propositions  obtained for the split quaternions as
in \cite{Ozdemir} to generalized them to similar results for the split octonions.
A split octonion is called \textit{spacelike, timelike or lightlike} if \ $%
N\left( a\right) <0,N\left( a\right) >0$ or $N\left( a\right) =0.$ If $%
N\left( a\right) =1,$ then $\ a$ is called \textit{the unit split octonion}.

\bigskip

\textbf{Spacelike octonions}

\bigskip

Let \ $a\in \mathbb{O}(1,1,-1)$ such that $N\left( a\right) =-1,$ be a
spacelike octonion. For the octonion $w=\frac{%
a_{1}f_{1}+a_{2}f_{2}+a_{3}f_{3}+a_{4}f_{4}+a_{5}f_{5}+a_{6}f_{6}+a_{7}f_{7}%
}{\sqrt{1+a_{0}^{2}}},$ we have $N\left( w\right) =$ $-1$ and $t\left(
w\right) =0,$ therefore $w^{2}=1.$ Denoting $\sinh \lambda =a_{0}$ and $%
\cosh \lambda =\sqrt{1+a_{0}^{2}},\lambda \in \mathbb{R},$ it results:

\begin{equation*}
a=e^{\lambda w}=\sinh \lambda +w\cosh \lambda .
\end{equation*}

If $a\in \mathbb{O}(1,1,-1)$ with $N\left( a\right) <0,$ we have \ $a=\sqrt{%
\left\vert N\left( a\right) \right\vert }(\sinh \lambda +w\cosh \lambda ).$

\textbf{Proposition 6.6. }\textit{We have that} $a^{n}=(\sqrt{\left\vert
N\left( a\right) \right\vert })^{n}(\sinh \lambda +w\cosh \lambda )$ \textit{%
for} $n$ \textit{odd and} $a^{n}=(\sqrt{\left\vert N\left( a\right)
\right\vert })^{n}(\cosh \lambda +w\sinh \lambda )$ \textit{for} $n$ \textit{%
even.} $\Box \medskip $

\textbf{Timelike octonions}

\bigskip

Let $a\in \mathbb{O}(1,1,-1)$ such that $N\left( a\right) =1,$ be a timelike
octonion. If $1-a_{0}^{2}>0,$ for the octonion $w=\frac{%
a_{1}f_{1}+a_{2}f_{2}+a_{3}f_{3}+a_{4}f_{4}+a_{5}f_{5}+a_{6}f_{6}+a_{7}f_{7}%
}{\sqrt{1-a_{0}^{2}}},$ we have $N\left( w\right) =$ $1$ and $t\left(
w\right) =0,$ therefore $w^{2}=-1$ Denoting $\cos \lambda =a_{0}$ and $\sin
\lambda =\sqrt{1-a_{0}^{2}},\lambda \in \mathbb{R},$ it results:

\textbf{Proposition 6.7.} \textit{With the above notations, we have the
Euler's formula:}%
\begin{equation*}
a=e^{\lambda w}=\cos \lambda +w\sin \lambda .
\end{equation*}

\textbf{Proof.} Indeed, $e^{\lambda w}=\overset{\infty }{\underset{i=1}{\sum
}}\frac{\left( \lambda w\right) ^{n}}{n!}$=$\overset{\infty }{\underset{i=1}{%
\sum }}\frac{\left( -1\right) ^{n}\lambda ^{2n}}{\left( 2n\right) !}+w%
\overset{\infty }{\underset{i=1}{\sum }}\frac{\left( -1\right) ^{n-1}\lambda
^{2n-1}}{\left( 2n-1\right) !}$=$\cos \lambda +w\sin \lambda .\Box $

If $a\in \mathbb{O}(1,1,-1)$ with $N\left( a\right) >0,$ it results \ $a=%
\sqrt{N\left( a\right) }\left( \cos \lambda +w\sin \lambda \right) .$

\textbf{Proposition 6.8.} \textit{We have that} $a^{n}=(\sqrt{N\left(
a\right) })^{n}(\cos n\lambda +w\sin n\lambda ).\Box $

\textbf{Proposition 6.9.}

\textit{1)If }$a\in \mathbb{O}(1,1,-1),$ \textit{it results} $a^{n}=(\sqrt{%
N\left( a\right) })^{n}(\cos n\lambda +w\sin n\lambda ).$

\textit{2)} \textit{The equation} $x^{n}=a$ \textit{has }$n$ \textit{roots:}
$\sqrt[2n]{\left\vert N\left( a\right) \right\vert }\left( \cosh \frac{%
\lambda }{n}+w\sinh \frac{\lambda }{n}\right) .\Box $

If $1-a_{0}^{2}<0,$ we have $w=\frac{%
a_{1}f_{1}+a_{2}f_{2}+a_{3}f_{3}+a_{4}f_{4}+a_{5}f_{5}+a_{6}f_{6}+a_{7}f_{7}%
}{\sqrt{a_{0}^{2}-1}},$ with $N\left( w\right) =$ $-1$ and $t\left( w\right)
=0,$ therefore $w^{2}=1.$ Denoting $\cosh \lambda =a_{0}$ and $\sinh \lambda
=\sqrt{a_{0}^{2}-1},\lambda \in \mathbb{R},$ it results:

\textbf{Proposition 6.10.} \textit{With the above notations, we have Euler's
formula:}%
\begin{equation*}
a=e^{\lambda w}=\cosh \lambda +w\sinh \lambda .
\end{equation*}

\textbf{Proof.} Indeed, $e^{\lambda w}=\overset{\infty }{\underset{n=0}{\sum
}}\frac{\left( \lambda w\right) ^{n}}{n!}$=$\overset{\infty }{\underset{n=0}{%
\sum }}\frac{\lambda ^{2n}}{\left( 2n\right) !}+w\overset{\infty }{\underset{%
n=0}{\sum }}\frac{\lambda ^{2n+1}}{\left( 2n+1\right) !}$=$\cosh \lambda
+w\sinh \lambda .\Box $

If $a\in a\in \mathbb{O}(1,1,-1)$ with $N\left( a\right) <0,$ it results \ $%
a=\sqrt{\left\vert N\left( a\right) \right\vert }(\cosh \lambda +w\sinh
\lambda ).$

\textbf{Proposition 6.11.}

\textit{1) If }$a\in \mathbb{O}(1,1,-1),$ \textit{then} $a^{n}=(\sqrt{%
\left\vert N\left( a\right) \right\vert })^{n}(\cosh n\lambda +w\sinh
n\lambda ).$

\textit{2)} \textit{The equation} $x^{n}=a$ \textit{has only one root:} $%
\sqrt[2n]{\left\vert N\left( a\right) \right\vert }\left( \cosh \frac{%
\lambda }{n}+w\sinh \frac{\lambda }{n}\right) .\Box \medskip $

\textbf{Remark 6. 12.} Using the above technique, De Moivre's formula and
Euler's formula can be easy proved for the octonion algebra $\mathbb{O}%
(\alpha ,\beta ,\gamma ),$ with $\alpha ,\beta ,\gamma \in \mathbb{R}%
\setminus \{0\}$ such that $\mathbb{O}(\alpha ,\beta ,\gamma )$ is split.
Thus, the study of algebraic equations in an arbitrary algebra $\mathbb{O}%
(\alpha ,\beta ,\gamma )$ with $\alpha ,\beta ,\gamma \in \mathbb{R}%
\setminus \{0\}$ was reduced to study of the corresponding algebraic
equation in one of the following two algebras: division octonion algebra $%
\mathbb{O}(1,1,1)$ or algebra $\mathbb{O}(1,1,-1)$.
\begin{equation*}
\end{equation*}

\begin{equation*}
\end{equation*}

\textbf{Conclusion.} In this paper, we used isomorphism between the real
quaternion algebras $\mathbb{H}(\gamma _{1},\gamma _{2})$ and $\mathbb{H}%
(1,1)$ or $\mathbb{H}(1,-1)$ to reduced the study of some algebraic
equations in an arbitrary algebra $\mathbb{H}\left( \gamma _{1},\gamma
_{2}\right) $ with $\gamma _{1},\gamma _{2}\in \mathbb{R}\setminus \{0\}$ to
study of the corresponding algebraic equation in one of the following two
algebras: division quaternion algebra or coquaternion algebra. The same
result was obtained for the generalized octonion algebra $\mathbb{O}(\alpha
,\beta ,\gamma )$. De Moivre's formula in generalized quaternion algebras
and generalized octonion division algebras was proved using this new method.

\vskip3mm %\newpage

Cristina FLAUT

{\small Faculty of Mathematics and Computer Science, Ovidius University,}

{\small Bd. Mamaia 124, 900527, CONSTANTA, ROMANIA}

{\small http://cristinaflaut.wikispaces.com/;
http://www.univ-ovidius.ro/math/}

{\small e-mail: cflaut@univ-ovidius.ro; cristina\_flaut@yahoo.com}

\medskip \medskip

Vitalii \ SHPAKIVSKYI

{\small Department of Complex Analysis and Potential Theory}

{\small \ Institute of Mathematics of the National Academy of Sciences of
Ukraine,}

{\small \ 3, Tereshchenkivs'ka st.,\ 01601 Kiev-4, UKRAINE}

{\small \ http://www.imath.kiev.ua/\symbol{126}complex/; \ e-mail:
shpakivskyi@mail.ru}


\begin{thebibliography}{99}
\bibitem{Carmody} K. Carmody, \textit{Circular and Hyperbolic Quaternions,
Octonions, and Sedenions--Further Results}, Appl. Math. Comput., \textbf{84}
(1997), 27-47.

\bibitem{Cockle} J. Cockle, \textit{On Systems of Algebra involving more
than one Imaginary}, Philosophical Magazine, \textbf{35} (3) (1849), 434-435.

\bibitem{Cho} E. Cho, \textit{De-Moivre's formula for quaternions}, Appl.
Math. Lett., \textbf{11} (6) (1998), 33-35.

\bibitem{Eilenberg} S. Eilenberg, I. Niven, \textit{The fundamental theorem
of algebra for quaternions}, Bull. Amer. Math. Soc., \textbf{50} (1944),
246-248.

\bibitem{Erdogdu} M. Erdo\u{g}du, M. \"{O}zdemir, \textit{Two-sided linear
split quaternionic equations with unknowns}, Linear and Multilinear Algebra,
DOI:10.1080/03081087.2013.851196

\bibitem{Flaut} C. Flaut, \textit{Some equations in algebras obtained by the
Cayley-Dickson process}, An. \c{S}t. Univ. Ovidius Constan\c{t}a, \textbf{9}%
(2)(2001), 45-68.

\bibitem{Flaut-St} C. Flaut, M. \c{S}tef\u{a}nescu, \textit{Some equations
over generalized quaternion and octonion division algebras}, Bull. Math.
Soc. Sci. Math. Roumanie, \textbf{52}(100), no. 4 (2009), 427-439.

\bibitem{Ivanov} S. Ivanov, S. Zamkovoy, \textit{Parahermitian and
paraquaternionic manifolds}, Diff. Geometry Appl., \textbf{23}, pp. 205-234.

\bibitem{Janovska} D. Janovsk\'{a}, G. Opfer, \textit{Linear equations and
the Kronecker product in coquaternions}, Mitt. Math. Ges. Hamburg \textbf{33}
(2013), 181-196.

\bibitem{Lam} T.Y. Lam, \textit{Quadratic forms over fields}, AMS,
Providence, Rhode Island, 2004.

\bibitem{Mamagani} A. B. Mamagani, M. Jafari, \textit{On Properties of
Generalized Quaternion Algebra}, J. Novel Appl. Sci., \textbf{2} (12)
(2013), 683-689.

\bibitem{Mierz} D. A. Mierzejewski, V. S. Szpakowski, \textit{On solutions
of some types of quaternionic quadratic equations}, Bull. Soc. Sci. Lett. L%
\'{o}d\'{z} \textbf{58}, Ser. Rech. D\'{e}form., \textbf{55} (2008), 49-58.

\bibitem{Mierz-1} D. A. Mierzejewski \textit{Linear manifolds in sets of
solutions of quaternionic polynomial equations of several types}, Adv. Appl.
Clifford Alg., \textbf{21} (2011), 417-428.

\bibitem{Mierz-2} D. Mierzejewski \textit{Spheres in sets of solutions of
quadratic quaternionic equations of some types}, Bull. Soc. Sci. Lett. L\'{o}%
d\'{z}, Ser. Rech. D\'{e}form., \textbf{60}(1)(2010), 49-58.

\bibitem{Niven} I. Niven, \textit{Equations in Quaternions}, Amer. Math.
Monthly, \textbf{48}(1941), 654-661.

\bibitem{Ozdemir} M. \"{O}zdemir, \textit{The roots of a split quaternion},
Appl. Math. Lett., \textbf{22} (2009) 258-263.

\bibitem{Pogoruy} A. Pogoruy, R. M. Rodrigues-Dagnino, \textit{Some
algebraic and analytical properties of coquaternion algebra}, Adv. Appl.
Clifford Alg., \textbf{20} (2010), 79-84.

\bibitem{Rosenfeld} B. A. Rosenfeld, \textit{A History of Non-Euclidean
Geometry}, Springer-Verlag, 1988.

\bibitem{Schafer} Schafer, R. D., \textit{An Introduction to Nonassociative
Algebras,} Academic Press, New-York, 1966.

\bibitem{Shpakivskyi} V. S. Shpakivskyi, \textit{Linear quaternionic
equations and their systems}, Adv. Appl. Clifford Alg., \textbf{21}(2011),
637-645.

\bibitem{Smith} W. D. Smith, \textit{Quaternions, octonions, and now,
16-ons, and 2n-ons; New kinds of numbers}, www. math. temple.edu/ {2dc}%
wds/homepage/nce2.ps, 2004.

\bibitem{Szpakowski} V. Szpakowski, \textit{Solution of general quadratic
quaternionic equations}, Bull. Soc. Sci. Lett. L\'{o}d\'{z}, Ser. Rech. D%
\'{e}form., \textbf{58}(2009), 45-58.

\bibitem{Tian} Y. Tian, \textit{Similarity and cosimilarity of elements in
the real Cayley-Dickson algebras}, Adv. Appl. Clifford Alg., \textbf{9} (1)
(1999), 61-76.

\medskip
\begin{equation*}
\end{equation*}
\end{thebibliography}
\end{document}